\documentclass[12pt, final]{amsproc}
\usepackage{amssymb}

\usepackage[cp850]{inputenc}
\usepackage{mathrsfs}
\usepackage{amscd}
\usepackage[notcite,notref]{showkeys}
\usepackage[all]{xy}
\usepackage{graphicx}
\usepackage{wrapfig}
\theoremstyle{plain}

\theoremstyle{remark}

\def\Lip{\operatorname{Lip}}

\def\L{\mathscr L}

\def\e{\varepsilon}

\def\R{\mathbb{R}}
\def\N{\mathbb{N}}
\def\H{\mathbb{H}}

\def\U{\mathscr{U}}
\def\V{\mathscr{V}}

\title{Fine structure of the homomorphisms of the lattice of uniformly continuous functions on the line}
\author{F\'elix Cabello S\'anchez}
\address{Departamento de Matem\'{a}ticas, UEx and IMUEx, 06071-Badajoz, Spain}
\email{fcabello@unex.es}

\thanks{Supported in part by DGICYT project MTM2016$\cdot$76958$\cdot$C2$\cdot$1$\cdot$P (Spain) and Junta de Extremadura programs GR$\cdot$15152 and  IB$\cdot$16056.}
\thanks{2010 Mathematics Subject Classification 46E05, 54C35.}

\begin{document}
\noindent{\footnotesize \boxed{\text{ Revised May 20, 2019. To appear in Positivity}}}\\[13pt]

\begin{abstract}
We provide a representation of the homomorphisms $U\longrightarrow \R$, where $U$ is the lattice of all uniformly continuous on the line. The resulting picture is sharp enough to 
describe the fine topological structure of the space of such homomorphisms.
\end{abstract}

\maketitle

\markboth{{\it F. Cabello S\'anchez}}{{\it Structure of the space of uniformly continuous functions on the line}}

\section{Introduction}
The purpose of this short note is to describe, as accurately as possible, the real-valued homomorphisms of the lattice of all uniformly continuous functions on the (half) line.

We denote by $U(X)$ the lattice of all real-valued uniformly continuous functions on $X$, which is invariably assumed to be a metric space. The sublattice of bounded functions is denoted by $U^*(X)$. When $X$ is the half-line $\H=[1,\infty)$  with the distance given by the absolute value we just write $U$ and $U^*$.

By a homomorphism of vector lattices we mean a linear map preserving joins and meets (equivalently, absolute values). Given a vector lattice $\mathscr L$,  we denote by $H(\mathscr L)$ the set of all homomorphisms $\phi:\mathscr L\longrightarrow\mathbb R$.

We are interested in $H(U)$ for two good reasons, apart from sheer curiosity.
The first and most obvious one is that, given an object of a category, the study of the homomorphisms against the ``simplest object''  in the category (if there is one) is interesting in its own right and often enlights the initial object. Quite clearly, $\R$ can be considered as the simplest vector lattice.

The study of the lattices of uniformly continuous functions and their homomorphisms has spurred a sustained, though moderate, interest for some time now; see the papers \cite{shi, ga-ja, javi, hus, hus-pul} and their references.

The second motivation springs from
the circle of ideas around the Samuel-Smirnov compactification, a classical construction in topology; see \cite{s} and \cite[Chapter 9, \S 41]{wi}. The space $H(\mathscr L)$ can be given the relative product topology it inherits from $\R^\mathscr L$.  The Samuel-Smirnov compactification of $X$ is then the  subspace of those homomorphisms $\phi: U^*(X)\longrightarrow\R$ that are unital in the sense that they send the function {\bf 1} to the number $1$.
This construction has attracted a considerable attention, 
even for very simple choices of the base space, such as  $X=\R^n$; see \cite{w, japos}

Very recently Garrido and Mero\~no \cite{g-m} have used the unital homomorphisms on $U(X)$ to construct a realcompactification of $X$ which plays the same role for general uniformly continuous functions than the Samuel-Smirnov compactification for those which are bounded.

Our modest contribution to this line of research is a description of the topological space $H(U)$ in the spirit of  Woods' \cite[Section 4]{w}. As we shall see, the most interesting homomorphisms are not unital and actually they vanish on every bounded function. So, somehow, there is a better life beyond 1.

We are aware of the fact that this is just one example and that most of the arguments presented here depend heavily on the peculiarities of the line.

This is compensated in part by the chief role played by the line amongst metric spaces as well as by the neat description of $H(U)$ that is achieved.

Moreover, as a byproduct, we compute the spaces of homomorphisms of other important lattices such as $U(\R)$ and $\Lip(\H)$ or $\Lip(\N)$.

\section{Elementary stuff}

This part contains a rather pedestrian description of $H(U)$, based on ultrafilters over the positive integers.
If $\mathscr L$ is a vector lattice, then $H(\mathscr L)$ is a subset of $\R^\L$.
This can be used to transfer to $H(\mathscr L)$ the product topology of $\R^\L$: a typical neighbourhood of $\phi$ has the form
$$
\{\psi\in H(\mathscr L): |\phi(f_i)-\psi(f_i)|<\varepsilon \text{ for $1\leq i\leq n$}\},
$$
where $f_1,\dots, f_n\in \mathscr L$ and $\varepsilon >0$. This is the only topology that we will consider on $H(\mathscr L)$.

\subsection{}
Let us begin with the observation that every $t\in\H$ gives rise to a unital homomorphism $\delta_t:U\longrightarrow\R$ by evaluation: $\delta_t(f)=f(t)$. 
It is \emph{really} easy to see that all unital homomorphisms arise as evaluation at some $t\in\H$.
 Indeed, put
$$
H_{\bf 1}=\{\phi\in H(U): \phi({\bf 1})=1\}.
$$
 Let us first check that these evaluations are dense in $H_{\bf 1}$. Pick $\phi\in H_{\bf 1}$, finitely many functions $f_i\in U$ and $\e>0$. We have to find a point $t\in \H$ such that
$$
|\phi(f_i)-f_i(t)|<\e\quad\quad(1\leq i\leq n).
$$
If we assume that no such $t$ exists, then, letting $c_i=\phi(f_i){\bf 1}$, we have
$$
\bigvee_{1\leq i\leq n}|c_i {\bf 1}- f_i|\geq \e{\bf 1}.
$$ 
But
$$
\phi\left(\bigvee_{1\leq i\leq n}|c_i {\bf 1}- f_i|\right)
=
\bigvee_{1\leq i\leq n}\phi(|c_i {\bf 1}- f_i|)
=
\bigvee_{1\leq i\leq n}|c_i \phi({\bf 1})- \phi(f_i)| =0,
$$ 
while $\phi(\e{\bf 1})=\e$, a contradiction.

Thus, given $\phi\in H_{\bf 1}$ we can find a net $(t_\alpha)$ such that $(\delta_{t_\alpha})$ converges to $\phi$ in $H(U)$. In particular we have
$$
\phi(f)=\lim_\alpha f(t_\alpha)
$$
for every $f\in U$. Taking $f={\frak t}$ as the identity of $\H$ and $t^*=\phi({\frak t})\geq 1$, we have
$t^*=\lim_\alpha t_\alpha
$
and so $\phi=\delta_{t^*}$. We have thus proved:
\medskip

\noindent$\bigstar$ \emph{ A homomorphism $\phi: U\longrightarrow \R$ has the form $\phi=c\delta_t$ for some $t\in\H$ and $0<c<\infty$ if and only if $\phi({\bf 1})>0.$ Otherwise $\phi$ vanishes on every bounded function.}
\medskip

The ``otherwise'' part is due to the fact that if $\phi({\bf 1})=0$, then $\phi(f)=0$ for every bounded $f\geq 0$ since $f\leq n{\bf 1}$ for some $n\in\N$ and so  $\phi(g)=0$ for every $g\in U^*$ since such a $g$ is the difference of two nonnegative functions in $U^*$.
\medskip

It is clear that the preceding proof relies on Heine-Borel theorem as it depends on the local compactness of the line. Let us remark, however, that $H(U^*)$ contains many unital homomorphisms which are not evaluations at points of $\H$. These form the Samuel-Smirnov compactification of $\H$; see \cite{s} and specially Section~4 in Woods' classical paper \cite{w}.

To see how these ``outer'' homomorphisms arise, take any uniformly separated sequence $(t_n)$, that is, such that $|t_n-t_k|>\e$ for some positive $\e$ and every $n\neq k$, and let $\mathscr U$ be a free ultrafilter on the integers. Then set
$$
\phi(f)=\lim_{\mathscr U(n)}f(t_n)\quad\quad(f\in U^*).
$$

The space $H(U)$ contains outer homomorphisms as well. These have to vanish on $U^*$ and, as we shall see, also at each function $f$ such that 
$t^{-1}f(t)\to 0$ as $t\to\infty$.

The main property of the half-line required here is that every $f\in U$ is Lipschitz for large distances: for every $\e>0$ there is a constant $L$, depending on $\e$ and $f$, such that
$$
|f(s)-f(t)|\leq L|s-t|\quad\text{provided}\quad |s-t|\geq\e;
$$
see \cite[Proposition~1.11]{b-l} or \cite[Lemma 2.2]{turkiyos} for the easy proof.
In particular the limit $L(f)=\limsup_{t\to\infty}t^{-1}|f(t)|$ is finite for every $f\in U$, which implies that for every $f\in U$ there is $c>0$ such that $c|f|\leq {\frak t}$ and so each homomorphism vanishing at $\frak t$ has to be zero. Thus, one can use $\phi({\frak t})$ to measure size in $H(U)$.

Note that there are (many) unbounded functions in $U$ such that 
$L(f)=0$, for instance $f(t)=t^\alpha$ for $0<\alpha<1$ or  $f(t)=\log t$.

Going back to $H(U)$, let $\U$ be a free ultrafilter on $\mathbb N$ and put
\begin{equation}\label{phiU}
\phi_\U(f)=\lim_{\mathscr U(n)}\frac{f(n)}{n}\quad\quad(f\in U).
\end{equation}
Clearly, $\phi_\U$ is correctly defined, belongs to $H(U)$, vanishes on $U^*$ and $\phi_\U(\frak t)=1$. Note that only the values of $f$ at the integers are used in the definition of $\phi_\U$. Now, consider the following subsets of $H(U)$:
\begin{align*}
H_\frak t&=\{\phi\in H(U):  \phi(\frak t)=1\},\\
H_\frak t^0&=\{\phi\in H(U): \phi(\frak t)=1 \text{ and }  \phi({\bf 1})=0\}.
\end{align*}
Every nonzero $\phi\in H(U)$ falls into $H_\frak t$ after renormalization: just take 
$\phi(\frak t)^{-1}\phi$.

\subsection{}\label{our}
Our immediate aim is to show that every $\phi\in H_\frak t^0$ comes from a free ultrafilter on $\N$, as in (\ref{phiU}).

Let us first check that the closure of the set $\{n^{-1}\delta_n:n\in\N\}$ in $H(U)$ contains $H_\frak t^0$.

This amounts to verifying that, given $\phi\in H_\frak t^0, f_i,\dots, f_k\in U$ and $\e>0$ there is $n\in\N$ such that
$$
\left|\phi(f_i)-\frac{f_i(n)}{n}\right|<\e\quad\quad(1\leq i\leq k).
$$
Assuming the contrary we have
$$
\bigvee_{1\leq i\leq k}\left|n\phi(f_i)-f_i(n)\right|\geq \e n\quad\quad(n\in\N).
$$
Letting $c_i=\phi(f_i)$ and taking into account that a uniformly continous function on $\H$ is bounded if and only if it bounded on $\N$,
we see that the function
$$
g= 0 \wedge \left(\bigvee_{1\leq i\leq k}\left|c_i{\frak t}-f_i\right|-\e {\frak t}\right)
$$
belongs to $U^*$, since it vanishes on $\mathbb N$, and satisfies
$$
\bigvee_{1\leq i\leq k}\left|c_i{\frak t}-f_i\right|\geq \e {\frak t}+g,
$$
which cannot be since
$$
\phi\left( \bigvee_{1\leq i\leq k}\left|c_i{\frak t}-f_i\right|\right)=
\bigvee_{1\leq i\leq k}\left|c_i\phi({\frak t})-\phi(f_i)\right|=0,
$$
while $\phi(\e {\frak t}+g)= \e\phi( {\frak t})+\phi(g)=\e$.

Now, let us fix $\phi\in H_\frak t^0$. It is clear that there is a filter on the integers, say $\mathscr F$, containing every set of the form
$
\{n\in \N: n^{-1}\delta_n\in V\},
$
where $V$ runs over the neighbourhoods of $\phi$ in $H(U)$. Now, if $\mathscr U$ is any ultrafilter refining $\mathscr F$, then $\phi=\phi_\U$ since for every $f\in U$ one has
$$
\phi(f)=\lim_{\mathscr F(n)}\frac{f(n)}{n}= \lim_{\mathscr U(n)}\frac{f(n)}{n}= \phi_\U(f).
$$
We therefore have:
\medskip

\noindent$\bigstar$ \emph{ Let $\phi\in H_\frak t$. If $\phi({\bf 1})>0$, then $\phi=t^{-1}\delta_t$ for some $t\in\H$. Otherwise there is a free ultrafilter $\U$ on the positive integers such that $\phi=\phi_\U$, as in $(\ref{phiU})$.}
\medskip

Hence, if $\phi$ vanishes at ${\bf 1}$ it also vanishes at every function with $L(f)=0$.

Of course we are proud of this statement.
However in its present form it cannot be used to detect when and why two ultrafilters induce the same homomorphism. This question leads to very interesting maths, as we will see in the next Section. 

\section{More advanced stuff}

In all what follows we denote by $\beta\N$ the Stone-\v{C}ech compactification of the positive integers and $\N^*=\beta\N\backslash\N$ will be the remainder. We understand each element of $\N^*$ as a free ultrafilter on $\N$ and each point of $\N$ as a fixed ultrafilter.
Let, as usual, $\ell_\infty$ denote the algebra of all bounded functions on $\mathbb N$, with the pointwise operations and order.

 As it is well-known, the ultrafilters on $\N$ are in exact correspondence with the algebra homomorphisms  $\ell_\infty\longrightarrow \R$ through the formula
$$
f\longmapsto \lim_{\U(n)}f(n).
$$ 

\subsection{}\label{declare}
Let us declare the ultrafilters $\U$ and $\mathscr V$ equivalent (and write $\U\approx\mathscr V$ for short) if they induce the same homomorphism on $U$, that is, when $\phi_\U=\phi_\mathscr V$. 
While two ultrafilters inducing the same homomorphism on $\ell_\infty$ actually agree, this is not the case for the notion of equivalence we have just introduced.

To see this, take $\U\in\N^*$ and put $\mathscr V=1+\U$, that is, the sets of $\mathscr V$ are obtained by translating those of $\U$ by a unit. Then for $f\in U$ we have
$$
\phi_\mathscr V(f)=\lim_{\mathscr V(n)}\frac{f(n)}{n}
=\lim_{\mathscr U(n)}\frac{f(n+1)}{n+1}=\lim_{\mathscr U(n)}\frac{f(n)}{n}= \phi_\mathscr U(f)
$$
since $f(n+1)-f(n)$ is bounded. Needless to say $\U$ and $\mathscr V$ are different as exactly one of them contains the set of even numbers.

Let us explain the notion of the image of an ultrafilter, which is implicit in the construction of the ``translate'' $1+\U$. Let $g:X\longrightarrow Y$ be a mapping, where $X$ and $Y$ are sets with no additional structure. If $\U$ is an ultrafilter on $X$, then the image of $\U$ under $g$ is the ultrafilter 
$$
\mathscr V= g[\U]=\{B\subset Y: g^\leftarrow[B]\in\U\}.
$$
Quite clearly, if $K$ is a compact Hausdorff space and $f:Y\longrightarrow K$ is any mapping, then 
one has
$$
\lim_{\mathscr V(y)}f(y)= \lim_{\mathscr U(x)}f(g(x)).
$$
In this way $1+\U$ is just the image of $\U$ under the translation $1+\bullet:\N\longrightarrow\N$ given by $(1+\bullet)(n)=1+n$.

Now, the idea is that if $\U\neq \mathscr V$
and $g:\N\longrightarrow\N$ increases fast enough, then $g[\U]\not\approx g[\V]$. 

Indeed, consider the function $2^\bullet:\N\longrightarrow\N$ defined by $2^\bullet(n)=2^n$. Let $\mathscr U$ and $\mathscr V$ be two different ultrafilters on $\N$.
We will prove that $\phi_{2^{\bullet}[\U]}\neq \phi_{2^{\bullet}[\V]}$, that is, that there is $f\in U$ such that
\begin{equation}\label{distinto}
\lim_{\mathscr U(n)}\frac{f(2^n)}{2^n}\neq  \lim_{\mathscr V(n)}\frac{f(2^n)}{2^n}
\end{equation}
Let $A$ be a witness set, so that $A$ belongs to $\U$ but not to $\V$. We define a Lipschitz $f:\H\longrightarrow\R$ as follows. First, we put
$$
f_0(2^n)=\begin{cases}
2^n& \text{if $n\in A$}\\
0 & \text{if $n\notin A$}
\end{cases}
$$
and $f_0(1)=1$ which corresponds to $n=0$.
Then we extend $f_0$ to a piecewise linear function on $\H$ thus: write $t\in [2^n,2^{n+1}]$ as $t=(1-s)2^n+s2^{n+1}$ with $0\leq s\leq 1$ and put
$$
f(t)=(1-s)f_0(2^n)+sf_0(2^{n+1}).
$$
The resulting function is Lipschitz (hence uniformly continuous) with Lipschitz constant at most
$$
\sup_{n\geq 0}\frac{|f(2^{n+1})-f(2^n)|}{|2^{n+1}-2^n|}\leq \sup_n\frac{2^{n+1}}{2^n}=2.
$$
Needless to say, for this $f$ the limit in the left-hand side of (\ref{distinto}) equals 1, while that on the right-hand side is 0 since $A^c$ belongs to $\mathscr V$. It's nice, isn't it?

\begin{figure}[h]\label{fig:all}
  \centering\vspace{-4pt}
  \includegraphics[width=154mm]{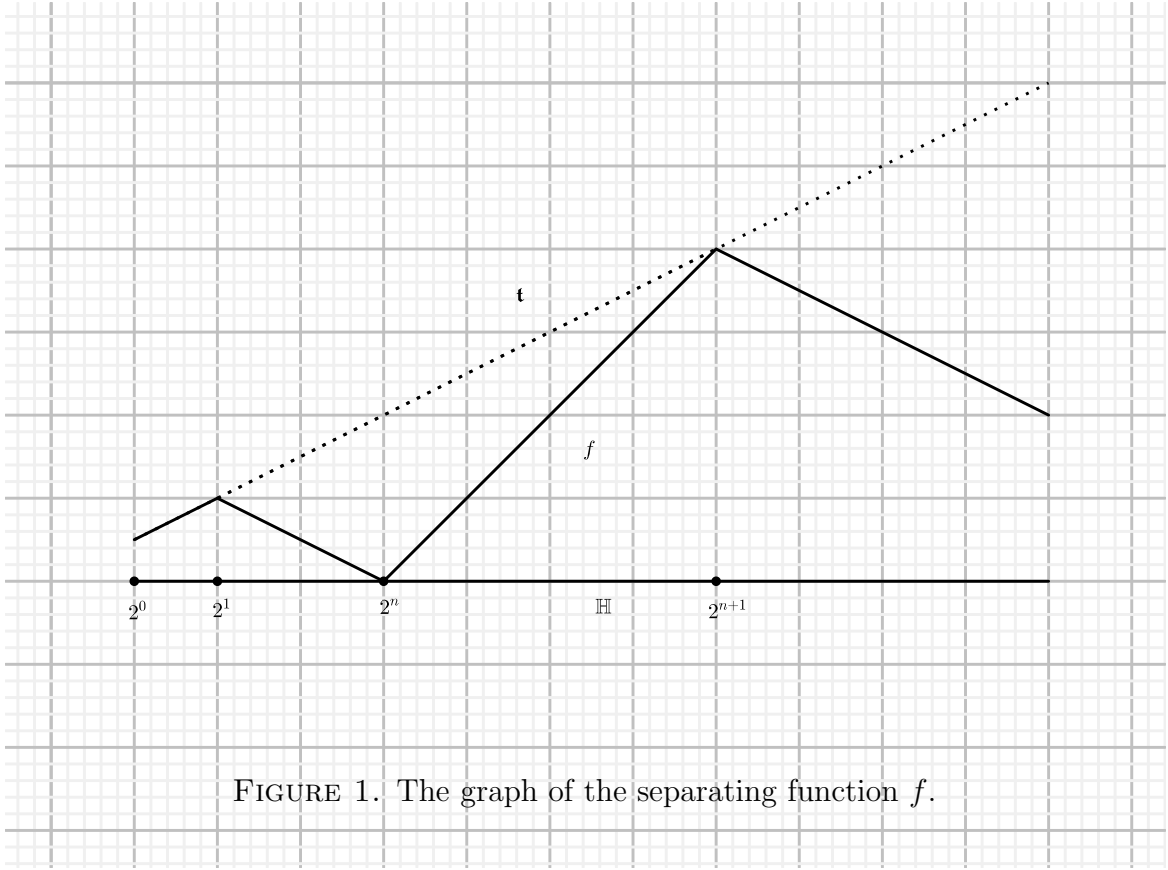}
  \vspace{-65pt}
  \caption[]
   {The graph of the separating function $f$.}
\end{figure}

\subsection{}\label{thepreceding}
The preceding observation is the key of the ensuing argument which allows us to give a neat description of the ``fiber'' $H_\frak t$ and so of $H(U)$. The fact that the exponential functions have exactly the growth-rate that is needed to separate ultrafilters is certainly a stroke of luck.

Let us denote by $\N_0$ the set of all nonnegative integers.
Every point $t\in[1,\infty)$ can be written as $t=c\cdot 2^n$ for some $n\in\N_0$ and $c\in[1,2]$. Let $m:[1,2]\times\N_0\longrightarrow[1,\infty)$ be the map sending $(c,n)$ to $c\cdot 2^n$  and notice that $m(c,n)=m(d,k)$ if and only if $c=2, d=1$ and $k=n+1$ or vice-versa. Composing $m$ with the mapping $[1,\infty)\longrightarrow H(U)$ sending $t$ to $t^{-1}\delta_t$ we obtain a map $\mu: [1,2]\times\N_0\longrightarrow H(U)$ namely $\mu(c,n)=m(c,n)^{-1}\delta_{m(c,n)}$, so
$$
\mu(c,n)(f)=\frac{f(c\cdot 2^n)}{c\cdot 2^n}\quad\quad(f\in U).
$$
Clearly, $\mu$ takes values in $
H_\frak t
$
which is a compact subset of $H(U)$ since it is homeomorphic to a closed subset of the product space
$
\prod_{f\leq\frak t}[0,1]
$. 
The argument appearing in {\bf \ref{our}} shows that the range of $\mu$ is dense in 
$
H_\frak t
$.
(Intermission: $
H_\frak t$ is a compactification of the half-line and also of the positive integers, with remainder $H^0_\frak t$ in both cases.)

Now we put Stone and \v{C}ech to work to obtain an extension $\tilde\mu: 
[1,2]\times\beta \N_0\longrightarrow H_\frak t$ so that, for every $c\in[1,2]$ and every $\U\in\N_0^*$ one has
$
\tilde{\mu}(c,\U)=\lim_{\U(n)}\mu(c,n)
$,
that is,
$$
\tilde{\mu}(c,\U)(f)= \lim_{\U(n)} \frac{f(c\cdot 2^n)}{c\cdot 2^n}\quad\quad(f\in U).
$$
Observe that the definition of $\tilde{\mu}$ guarantees continuity in the second variable, but not joint continuity that we now check ``by hand''.

\subsection{}\label{joint}
Recall that the topology of $\beta\N_0$ comes from $\ell_\infty(\N_0)$ in the sense that, given $\U\in\beta\N_0$, the sets of the form
$$
\{\mathscr V: |f(\mathscr V)-f(\U)|\leq\e\}\quad\quad(f\in\ell_\infty,\e>0),
$$
where $g(\mathscr W)$ is interpreted as the limit of $g(n)$ along $\mathscr W$, form a subbase of the topology at $\U$.

Thus, in order to stablish the continuity of $\tilde{\mu}$ at $(c,\U)$, it suffices to prove that, given $f\in U$ and $\e>0$, there exist $g\in\ell_\infty(\N_0)$ and $\delta>0$ so that $|c-d|<\delta$ and $|g(\mathscr V)-g(\U)|<\delta$ imply
$$
\left|\lim_{\U(n)} \frac{f(c\cdot 2^n)}{c\cdot 2^n}- 
\lim_{\mathscr V(n)} \frac{f(d\cdot 2^n)}{d \cdot 2^n}
\right|<\e.
$$
Let us assume that $\U$ is free. The other case is easier.
As the reader may guess we take $g(n)=f(c2^n)/(c2^n)$. Also, let $L$ be such that
$|f(s)-f(t)|\leq L|s-t|$ provided $|s-t|\geq 1$ and suppose $|d-c|<\e$ and
 $|g(\mathscr V)-g(\U)|<\e$. Then
 $$
\left|\lim_{\U(n)} \frac{f(c2^n)}{c 2^n}- 
\lim_{\mathscr V(n)} \frac{f(d 2^n)}{d 2^n}
\right|\leq 
\underbrace{\left|\lim_{\U(n)} \frac{f(c2^n)}{c 2^n}- 
\lim_{\mathscr V(n)} \frac{f(c 2^n)}{d  2^n}
\right|}_{(\dagger)}
+
\underbrace{\left|\lim_{\mathscr V(n)} \frac{f(c2^n)}{d  2^n}-
\lim_{\mathscr V(n)} \frac{f(d 2^n)}{d  2^n}
\right|}_{(\ddagger)}\!.
$$
Now,
$$
(\dagger)= \left|g(\U)-\frac{c}{d}g(\mathscr V)\right|\leq 
 \left|g(\U)-g(\mathscr V)\right|+  \left|g(\mathscr V)-\frac{c}{d}g(\mathscr V)\right|
 \leq \e+ \e g(\mathscr V)\leq \e(1+L),
$$

$$
(\ddagger)\leq \lim_{\mathscr V(n)}  \left| \frac{f(c2^n)-f(d2^n)}{d  2^n} \right|
\leq \lim_{\mathscr V(n)} \frac{L|c-d|2^n}{d2^n}\leq L\e,
$$
and so $\tilde{\mu}$ is a continuous mapping \emph{onto} $H_\frak t$.

\subsection{}\label{identify}
We have just seen that $\tilde\mu: 
[1,2]\times\beta \N_0\longrightarrow H_\frak t$ is a continuous surjection. The next task is to find out when two points of $[1,2]\times\beta \N_0$ have the same image in $H_\frak t$. Most of the work has been already done in {\bf \ref{declare}}, though in an implicit way.

Pick $(c,\U)$ and $(d,\mathscr V)$ in $[1,2]\times\beta \N_0$.

 We claim that  $\tilde{\mu}(c,\U)=\tilde{\mu}(d,\mathscr V)$ if and only if 
 $c=2, d=1$ and $\mathscr V=1+\U$ or vice-versa. The meaning of this equality was explained in {\bf \ref{declare}}.
 
 Let us first show that $\tilde{\mu}(2,\U)=\tilde{\mu}(1,1+\mathscr U)$. Take $f\in U$. Then
 $$ 
\tilde{\mu}(2,\U)(f)= 
\lim_{\U(n)} \frac{f(2\cdot 2^n)}{2\cdot 2^n}=
\lim_{\U(n)} \frac{f(2^{1+n})}{2^{1+n}}=
\lim_{(1+\U)(n)} \frac{f(1\cdot 2^{n})}{1\cdot 2^{n}}=
\tilde{\mu}(1,1+\mathscr U)(f),
 $$
as required. To check the converse we may assume $1\leq c\leq  d<2$ since otherwise we could replace $(2,\U)$ by $(1,1+\U)$ and/or $(2,\V)$ by $(1,1+\V)$.
Let us consider the case where $\U\neq \V$. Then there is $A\in \U$ which does not belong to $\V$ and so $A^c\in \V$. Here $A^c=\mathbb N_0\backslash A$ is the complement of $A$. Write
$$
\mathbb H=\bigsqcup_{n\geq 0}[2^n, 2^{n+1}).
$$
Each interval $[2^n, 2^{n+1})$ contains exactly one point of the form $c2^n$ with $n\in\N_0$ and another one of the form $d2^n$.
We define an increasing sequence $(p_n)_{n\geq 0}$ taking
$$
p_n=\begin{cases}
c2^n &\text{if $n\in A$}\\
d2^n &\text{if $n\notin A$}
\end{cases}
$$
Note that
$$
|p_{n+1}-p_n|\geq \min(c,d)2^{n+1}- \max(c,d)2^{n} =
c2^{n+1}- d2^{n} = (2c-d)2^n.
$$
Since $2c-d>0$ there exists a Lipschitz (hence uniformly continuous) $f:\mathbb H\longrightarrow\R$ such that
$$
f(t)=\begin{cases}
c2^n & \text{if $t=c2^n$ and $n\in A$}\\
0 & \text{if $t=d2^n$ and $n\in A^c$}
\end{cases}
$$
namely the function whose graph is the polygonal joining $(p_0,q_0), (p_1,q_1); (p_2,q_2)...$ with $q_n=p_n$ if $n\in A$ and $q_n=0$ otherwise. For this $f$ one clearly has $\tilde{\mu}(c,\U)(f)=1$, while $\tilde{\mu}(d,\mathscr V)(f)=0$, so
$\tilde{\mu}(c,\U)\neq\tilde{\mu}(d,\mathscr V)$.

Finally, if $\U=\V$, but $c\neq d$, then one easily finds a Lipschitz $f$ such that 
$$
f(t)=\begin{cases}
c2^n & \text{if $t=c2^n$}\\
0 & \text{if $t=d2^n$}
\end{cases}
$$
for every $n$ from where it follows that $\tilde{\mu}(c,\U)\neq\tilde{\mu}(d,\mathscr U)$. Since the case where $\U$ or $\V$ are fixed is trivial this leads to the following description of $H_{\frak t}$. Note that every continuous surjection between Hausdorff compacta is automatically a quotient map (cf. Willard \cite[Chapter~3, \S 9]{w}).


\medskip

\noindent$\bigstar$
\emph{The fiber $H_{\frak t}$ is homeomorphic to the quotient obtained from $[1,2]\times\beta\N_0$ after identifying each point of the form $(2,\U)$ with $(1,1+\U)$. The map sending the class of $(c,\U)$ to the homomorphism defined by the formula
$$
\phi(f)=\lim_{\U(n)} \frac{f(c\cdot 2^n)}{c\cdot 2^n}\quad\quad(f\in U)
$$
is a homeomorphism.}
\medskip

To complete our picture of $H(U)$, note that $H(U)\backslash\{0\}$ is homeomorphic to $H_\frak t\times(0,\infty)$: the map $(\phi,\lambda)\longmapsto \lambda\phi$ is continuous, with continuous inverse given by $\varphi\longmapsto(\varphi(\frak t)^{-1}\varphi, \varphi(\frak t))$.

On the other hand since for each $f\in U$ there exist $c>0$ such that $c|f|\leq\frak t$ we see that the sets
$
\{\phi\in H(U): \phi(\frak t)<\e\}
$
form a base of neighbourhoods of $0$. Hence:

\medskip

\noindent$\bigstar$
\emph{The space $H(U)$ is homeomorphic to the quotient of
$[1,2]\times\beta\N_0\times(0,\infty)$ with one point ${\bf 0}$ added, where we identify points of the form $(2,\U,\lambda)$ and $(1,1+\U,\lambda)$ and the  
neighbourhoods of the point ${\bf 0}$ are those sets containing a subset of the form $\{(c,\U,\lambda): \lambda<\e \}$ for some $\e>0$ together with the point ${\bf 0}$.}

\subsection{}\label{other}
It is clear that everything what has been said about $U$ applies verbatim to $\Lip(\mathbb H)$, the lattice of Lipschitz functions on the half-line. Hence
the  spaces of homomorphisms of $\Lip(\mathbb H)$ and $U$ agree, in the sense that each homomorphism $\Lip(\mathbb H)\longrightarrow\R$ is the restriction of a unique $\phi\in H(U)$. We refer the reader to Chapter 5 of Weaver booklet \cite{weaver} for basic information about Lipschitz lattices.

Also, since the line $\R$ can be obtained by ``gluing'' two half-lines,  $H(U(\R))$ can be easily computed using two copies of $H(U)$ and the same applies to $\Lip(\R)$. We will not give the details.

\medskip

Finally, let us describe the homomorphisms $\Lip(\N)\longrightarrow\R$, where $\N$ carries the metric inherited from $\R$. Let $E:\Lip(\N)\longrightarrow U$ be the linear map sending $f$ into the piecewise linear function on $\H$ that interpolates $f$ on $\N$. If we consider $\ell_\infty$ as the set of bounded functions in $\Lip\N$, then $E$ maps $\ell_\infty$ to $U^*$ and we cannot help to display the following commutative diagram of linear maps
$$
\xymatrix{
0\ar[r] & U^* \ar[r] & U \ar[r] & U/U^* \ar[r] & 0\\ 
0\ar[r] & \ell_\infty \ar[r]  \ar[u]^E & \Lip(\N) \ar[r] \ar[u]^E & \Lip(\N)/\ell_\infty \ar[r] \ar@{=}[u] & 0\\ 
}
$$
Here, the rows are exact and the equal sign on the right reflects the fact that every $f\in U$ agrees with one of the form $E(g)$ for some $g\in \Lip\N$ modulo a bounded function: actually one can take $g=E(f|_\N)$.
Now, let $L=H(\Lip \N)$ and
\begin{itemize}
\item $L_{\frak t}=\{\phi\in L: \phi({\frak t})=1\}$,
\item $L_{\frak t}^0=\{\phi\in L: \phi({\frak t})=1 \text{ and } \phi({\bf 1})=0\}$.
\end{itemize}
If $\phi\in L_{\frak t}$, then either $\phi({\bf 1})>0$, in which case $\phi=n^{-1}\delta_n$ for some integer $n$, or $\phi$ vanishes on every bounded function and so it factors throught the quotient $\Lip(\N)/\ell_\infty=U/U^*$. If so, there is $\U\in \N^*$ and $c\in[1,2]$ such that
$$
\phi(f)=\lim_{\U(n)}\frac{Ef(c2^n)}{c2^n}= \lim_{\U(n)}\frac{Ef([c2^n])}{[c2^n]},
$$
where $[\cdot]$ is the integer part function. Hence $L_{\frak t}^0$ is homeomorphic to $H_{\frak t}^0$.

\section{Coda}
This note lived a hard life until the referee it was looking for came along. 
In the meantime, we have explored in \cite{quiz} the consequences of the research reported here, encountered nice descriptions of the Samuel-Smirnov compactification of the line in \cite[Lemma 2.1]{flows} and \cite[Theorem 2.1]{turkiyos}, and registered some curious connections between the homomorphisms on $\Lip(\N)$ and the so-called density measures (cf. \cite[Theorem 3.1]{density}).


\end{document}